\documentclass[11pt]{article}
\usepackage{amsfonts}
\usepackage{amsfonts}
\usepackage{indentfirst,latexsym,bm,amsmath,amsthm}
\usepackage{amssymb}

\newtheorem{remark}{Remark}[section]
\newtheorem{lemma}{Lemma}[section]

\newtheorem{theorem}{Theorem}[section]
\newtheorem{prop}{Proposition}[section]

\newtheorem{definition}{Defintion}[section]

\begin{document}
\title{\bf Large deviations for invariant measures of  SPDEs with two reflecting walls}
\author{Tusheng Zhang\\
School of Mathematics, University of Manchester\\
Oxford Road, Manchester M13 9PL,\\
 England, U.K. \\
 Email: tusheng.zhang@manchester.ac.uk}
\date{}
\maketitle

\begin{abstract}
In this article, we established a large deviation principle for  invariant measures of solutions of   stochastic partial differential equations with two reflecting  walls
driven by space-time white noise. \\


\end{abstract}

\section{Introduction}
\setcounter{equation}{0}
Consider reflected stochastic partial differential equations (SPDEs) of the following type:
\begin{eqnarray}\label{I1}
&&\frac{\partial{u_{\varepsilon}(x,t)}}{\partial{t}}=\frac{\partial^2{u_{\varepsilon}(x,t)}}{\partial{x^2}}-\alpha u_{\varepsilon}(x,t)+f\big(x,u_{\varepsilon}(x,t)\big)
\nonumber\\
&&\quad\quad\quad \quad\quad +\varepsilon \sigma\big(x,u_{\varepsilon}(x,t)\big)\dot{W}(x,t)+\eta_{\varepsilon}(x,t)-\xi_{\varepsilon}(x,t)\\
&& K_1(x)\leq {u^{\varepsilon}}(x,t)\leq K_2(x)
\end{eqnarray}
in $(x,t)\in{Q}:=[0,1]\times\mathbb{R}_+$ while
$K_1(x)\leq {u^{\varepsilon}}(x,t)\leq K_2(x)$. Here $\dot{W}$ is a space-time white noise.  When $u_{\varepsilon}(x,t)$
hits $K_1(x)$ or $K_2(x)$, the additional forces are added to
prevent $u_{\varepsilon}$ from leaving $[K_1,\ K_2]$. These forces are
expressed by random measures $\xi_{\varepsilon}$ and $\eta_{\varepsilon}$ in equation (\ref{I1}) which play
a similar role as the local time in the usual Skorokhod equation constructing Brownian motions with reflecting barriers.
\vskip 0.3cm
Parabolic SPDEs with reflection are natural extension of the widely studied deterministic parabolic
obstacle problems. They also can be used to model fluctuations of an interface near a wall, see Funaki and
Olla \cite{FO}.
In recent years, there is a growing interest on the study of SPDEs with reflection. Several works are
devoted to the existence and uniqueness of the solutions. In the case
of a constant diffusion coefficient and a single reflecting barrier $K_1=0$,
Naulart and Pardoux \cite{NP} proved the existence and uniqueness of the solutions.
In the case of a non-constant diffusion coefficient and a single reflecting barrier $K_1=0$,
the existence of a minimal solution was obtained by
Donati-Martin and Pardoux \cite{DP1}. The existence and particularly the uniqueness of the solutions for
a fully non-linear SPDE with reflecting barrier $0$ were solved by Xu and Zhang \cite{XZ}.
In the case of double reflecting barriers, Otobe \cite{OS} obtained the
existence and uniqueness of the solutions of a SPDE
driven by an additive white noise.
\vskip 0.3cm
In addition to the existence and uniqueness, various other
properties of the solution have been studied by several authors, see
Donati-Martin and Pardoux \cite{DP2}, Zambotti \cite{ZA}, Dalang et al. \cite{DMZ} and Zhang \cite{ZW}.
\vskip 0.3cm
The purpose of this paper is to establish a large deviation principle for invariant measures of  the solutions of
fully non-linear SPDEs  with two reflecting walls (\ref{I1}). Large deviations for  invariant measures of  the solutions of  SPDEs were previously studied in \cite{S} and \cite{CR}. Our approach will be  along the same lines as that in \cite{S} and \cite{CR}. However, the extension is non-trivial. The extra difficulty arises from the appearance of the random measures ( local times) $\eta_{\varepsilon}$ and  $\xi_{\varepsilon}$ in the equation (\ref{I1}). We need to carefully analyze the local time terms in the skeleton equations and provide some uniform  estimates for the penalized   approximating equations.
\vskip 0.3cm
The rest of the paper is organized as follows. In Section 2 we introduce the SPDEs with reflecting walls and state the precise conditions on the coefficients. In Section 3 we recall some results on the deterministic obstacle problems which will be used later. In Section 4, we study the skeleton equations and  the rate functional. We  provide some estimates
for the extra measures (local times) in the equation and prove  equivalent characterizations of the rate functional.
In Section 5, we prove the exponential tightness for the invariant measures. The main result is stated in Section 6. The lower bound of the large deviation is established in Section 7 and upper bound is obtained in Section 8.

\section{Reflected SPDEs}
\setcounter{equation}{0}
In this section, we introduce reflected stochastic partial differential equations (SPDEs)  and state the precise conditions on the coefficients.

Consider the following SPDE with two reflecting walls:
\begin{equation}\label{0.1}
\left \{\begin{array}{ll}
        \frac{\partial{u}}{\partial{t}}=\frac{\partial^2{u}}{\partial{x^2}}-\alpha u+f(x,u(x,t))+\sigma(x,u(x,t))\dot{W}(x,t)+\eta-\xi;\\
        \frac{\partial u}{\partial x}(0,t)=0,\ \frac{\partial u}{\partial x}(1,t)=0,\ \ {\rm for}\ t\geq0;\\
        u(x,0)=u_0(x)\in{C}([0,1]); K_1(x)\leq{u}_0(x)\leq{K_2}(x), \\
        K_1(x)\leq{u}(x,t)\leq{K_2}(x),\ \ {\rm for}\ (x,t)\in{Q},
\end{array}
\right.
\end{equation}
here $W(x,t)$ is a space-time Brownian sheet on a filtered probability space
$(\Omega,P,\mathcal {F};\mathcal {F}_t, t\geq 0 )$.
\vskip 0.3cm
Throughout the paper, the reflecting walls $K_1(x),\ K_2(x)$ are assumed to be continuous
functions satisfying \\
(H1) $K_1(x)<0<K_2(x)$ for $x\in [0,1]$; \\
(H2)
$\frac{\partial^2{K_i}}{\partial{x^2}}\in
L^2([0,1])$, where
$\frac{\partial^2{}}{\partial{x^2}}$ are interpreted in a distributional sense.\\
 \\
 Introduce the following conditions: \\
(F1) the coefficients: $f, \sigma:
[0,1]\times\mathbb{R}\rightarrow\mathbb{R}$
are bounded and there exists $C>0$ such that
$$|{f}(x,y)-f(x,\hat{y})|+|\sigma(x,y)-\sigma(x,\hat{y})|\leq{C}|y-\hat{y}|,$$
for $x\in[0,1]$ and $y, \hat{y}\in \mathbb{R}$.\\
(F2) $\sigma (x,y)$ is continuous in both variables and there exists $m>0$ such that $|\sigma (x,y)|\geq m$. \\
\begin{remark}
Here in the equation (\ref{0.1}) we choose the Neumann boundary condition for the Laplacian operator for convenience.
The results in this paper are also valid for other boundary conditions, e.g, periodic boundary condition, Dirichlet boundary condition.
\end{remark}
The following is the definition of the solution of a SPDE with two reflecting walls $K_1,\ K_2$.
\begin{definition} A triplet
$(u,\eta,\xi)$  is a solution to the
SPDE
(\ref{0.1}) if\\
(i) $u=\{u(x,t);(x,t)\in{Q}\}$ is a continuous, adapted random field
(i.e., $u(x,t)$ is $\mathcal {F}_t$-measurable $\forall$ $t\geq0,
x\in[0,1]$) satisfying $K_1(x)\leq{u}(x,t)\leq{K_2}(x)$, a.s;\\
(ii) $\eta(dx,dt)$ and $\xi(dx,dt)$ are positive and adapted (i.e.
$\eta(B)$ and $\xi(B)$ is $\mathcal {F}_t$-measurable if
$B\subset[0,1]\times[0,t]$) random measures on $[0,1]\times\mathbb{R}_+$
satisfying
$$\eta\big([0,1]\times[0,T]\big)<\infty,\ \xi\big([0,1]\times[0,T]\big)<\infty$$
for  $T>0$;\\
(iii) for all $t\geq0$ and $\phi\in{C^\infty}[0,1]$(the set of smooth
functions) with $\frac{\partial \phi}{\partial x}(0)= \frac{\partial \phi}{\partial x}(1)=0$ we
have
\begin{eqnarray}
&&\big(u(t),\phi\big)-\int_0^t(u(s),\phi^{''})ds-\int_0^t\big(f(\cdot,u(s)),\phi\big)ds-
\int_0^t\int_0^1\phi(x)\sigma(x,u(x,s))W(dx,ds)\nonumber\\
&=&\big(u_0,\phi\big)-\alpha \int_0^t(u(s),\phi)ds+\int_0^t\int_0^1\phi(x)\eta(dx,ds)-\int_0^t\int_0^1\phi(x)\xi(dx,ds),\ a.s,
\end{eqnarray}
where $(,)$ denotes the inner product in $L^2([0,1])$ and $u(t)$ denotes $u(\cdot,t)$;\\
(iv)$$\int_Q\big(u(x,t)-K_1(x)\big)\eta(dx,dt)=\int_Q\big(K_2(x)-u(x,t)\big)\xi(dx,dt)=0.$$
\end{definition}

\section{Deterministic obstacle problems}
\setcounter{equation}{0}
Let $K_1,\ K_2$ be as in Section 2 and $u_0\in{C}([0,1])$ with $K_1(x)\leq{}u_0(x)\leq{}K_2(x)$.
Let $v(x,t)\in{}C(Q)$ with $v(x,0)=u_0(x)$.
Consider a deterministic PDE with two reflecting walls:
\begin{equation}\label{1.1}
\left \{\begin{array}{ll}
       \frac{\partial{z(x,t)}}{\partial{t}}-\frac{\partial^2{z(x,t)}}{\partial{x^2}}+\alpha z(x,t)=\eta(x,t)-\xi(x,t);\\
        \frac{\partial z}{\partial x}(0,t)=\frac{\partial z}{\partial x}(1,t)=0,\ \ {\rm for}\ t\geq0;\\
        z(x,0)=0,\ \ {\rm for}\ x\in[0,1];\\
        K_1(x)-v(x,t)\leq{z}(x,t)\leq{K_2}(x)-v(x,t),\ \ {\rm for}\ (x,t)\in{Q}.
\end{array}
\right.
\end{equation}
We first present a precise definition of the solution for equation (\ref{1.1}).
\begin{definition} A triplet $(z,\eta,\xi)$ is called a
solution to the PDE (\ref{1.1}) if\\
(i) $z={z(x,t);(x,t)\in{Q}}$ is a continuous function satisfying
$K_1(x)\leq{z}(x,t)+v(x,t)\leq{K_2}(x),\ z(x,0)=0$;\\
(ii) $\eta(dx,dt)$ and $\xi(dx,dt)$ are measures on
$[0,1]\times\mathbb{R}_+$ satisfying
$$\eta\big([0,1]\times[0,T]\big)<\infty,\  \xi\big([0,1]\times[0,T]\big)<\infty$$
for $T>0$;\\
(iii) for all $t\geq0$ and $\phi\in{C^\infty}[0,1]$ with $\frac{\partial \phi}{\partial x}(0)=\frac{\partial \phi}{\partial x}(1)=0$ we have
\begin{eqnarray}\label{1.2}
     &&\big(z(t),\phi\big)-\int_0^t(z(s),\phi^{''})ds +\alpha \int_0^t (z(s), \phi )ds\nonumber\\
     &=&\int_0^t\int_0^1\phi(x)\eta(dx,ds)-\int_0^t\int_0^1\phi(x)\xi(dx,ds),
\end{eqnarray}
where $z(t)$ denotes $z(\cdot,t)$.\\
(iv)
\begin{eqnarray*}
     &&\int_Q\big(z(x,t)+v(x,t)-K_1(x)\big)\eta(dx,dt) \\
          &=&\int_Q\big(K_2(x)-z(x,t)-v(x,t)\big)\xi(dx,dt)\\
          &=&0.
\end{eqnarray*}
\end{definition}

The following result is the existence and uniqueness of the solutions.We refer the reader to \cite{YZ1} for the proof.
\begin{theorem} Suppose (F1) holds. Equation (\ref{1.1}) admits a unique solution
$(z,\eta,\xi)$.
\end{theorem}
\hfill$\Box$
\begin{remark} Let $\hat{z}$ be the solution to equation (\ref{1.1}) replacing $v$
by $\hat{v}$, where $\hat{v}(x,t)$ is another continuous function on
Q. It is proved in Otobe \cite{OS} that
\begin{eqnarray}\label{1.3}
\|z-\hat{z}\|_\infty^T\leq\|v-\hat{v}\|_\infty^T,
\end{eqnarray}
for any $T>0$, where
$\|\omega\|_\infty^T:=\underset{x\in[0,1],t\in[0,T]}{\sup}|\omega(x,t)|$.
\end{remark}

\section{Skeleton equations and rate functional}
\setcounter{equation}{0}
The Cameron-Martin space associated with
the Brownian sheet $\{W(x,t), x\in[0,1], t\in \mathbb{R}_+\}$ is
given by
$$\mathcal{H}=\{h=\int_0^\cdot\int_0^\cdot \dot{h}(x,s)dxds; \int_0^T\int_0^1 \dot{h}^2(x,s)dxds<\infty, T>0\}.$$
For $h=\int_0^\cdot\int_0^\cdot \dot{h}(x,s)dxds\in \mathcal{H},$
consider the following reflected deterministic PDE (the skeleton
equation):
\begin{eqnarray}\label{3.1}\left\{
\begin{array}{l}
\frac{\partial u^h(x,t)}{\partial t}-\frac{\partial^2 u^h(x,t)}{\partial x^2}+\alpha u^h(x,t)=f(x,u^h(x,t))+\sigma(x,u^h(x,t))\dot{h}(x,t)+\eta^h-\xi^h;\\
K_1(x)\leq u^h(x,t)\leq K_2(x);\\
\int_0^T\int_0^1 (u^h(x,t)-K_1(x))\eta^h(dx,dt)=\int_0^T\int_0^1 (K_2(x)-u^h(x,t))\xi^h(dx,dt)=0, T>0\\
u^h(\cdot,0)=u_0;\\ \frac{\partial u^h}{\partial x}(0,t)=\frac{\partial u^h}{\partial x}(1,t)=0.
\end{array}\right.
\end{eqnarray}
\begin{theorem}
Assume (F1). Then, equation (\ref{3.1}) admits a unique solution. Moreover, the measures $\eta^h,\xi^h$ are absolutely continuous w.r.t. the Lebesgue measure $dx\times dt$ and for $T>0$,
\begin{equation}\label{3.2}
\int_0^T\int_0^1 |\dot{\eta}^h(x,t)|^2dxdt\leq C_{1,T,h}, \quad \int_0^T\int_0^1 |\dot{\xi}^h(x,t)|^2dxdt\leq C_{2,T,h},
\end{equation}
where $C_{1,T,h}, C_{2,T,h}$ are constants only depending on the bounds of $f$, $\sigma$, and the norm $||h||$. Furthermore, if $||h||_2^{\infty}=\int_0^{\infty}\int_0^1 \dot{h}^2(x,s)dxds<\infty $, then we have
\begin{equation}\label{3.2-1}
\int_0^Te^{-\alpha (T-t)}\int_0^1 |\dot{\eta}^h(x,t)|^2dxdt\leq C(1+||h||_2^{\infty}),
\end{equation}
\begin{equation}\label{3.2-2}
 \quad \int_0^Te^{-\alpha (T-t)}\int_0^1 |\dot{\xi}^h(x,t)|^2dxdt\leq C(1+||h||_2^{\infty}),
\end{equation}
where $C$ is a constant independent of $T$.
\end{theorem}
\noindent\textbf{Proof.}\quad
Consider the approximating equations:
\begin{eqnarray}\label{3.3}
       \frac{\partial{u_{\varepsilon,\delta}^h(x,t)}}{\partial{t}}&=&\frac{\partial^2{u_{\varepsilon,\delta}^h(x,t)}}{\partial{x^2}}
       -\alpha u_{\varepsilon,\delta}^h(x,t)+f(u_{\varepsilon,\delta}^h(x,t))
       +\sigma(u_{\varepsilon,\delta}^h(x,t))\dot{h}(x,t)\nonumber\\
       && +\frac{1}{\delta}(u_{\varepsilon,\delta}^h(x,t)-K_1(x))^- -\frac{1}{\varepsilon}(u_{\varepsilon,\delta}^h(x,t)-K_2(x))^+; \\
       u_{\varepsilon,\delta}^h(x,0)&=&u_0(x); \nonumber \\
       u_{\varepsilon,\delta}^h(0,t)&=&u_{\varepsilon,\delta}^h(1,t)=0.\nonumber
\end{eqnarray}
Here for simplicity, we write $f(u_{\varepsilon,\delta}^h(x,t))$ for $f(x,u_{\varepsilon,\delta}^h(x,t))$ and
       $\sigma(u_{\varepsilon,\delta}^h(x,t))$ for $\sigma(x,u_{\varepsilon,\delta}^h(x,t))$.
As shown in the SPDE case in \cite{OS},\cite{YZ1},  the solution $u_{\varepsilon,\delta}^h$ of equation (\ref{3.3}) converges to the  unique solution $u^h$ of (\ref{3.1}) as $\delta, \varepsilon \rightarrow 0$. Moreover
\begin{equation}\label{3.4}
\eta^h(dx,dt)=\lim_{\varepsilon\rightarrow 0}\lim_{\delta\rightarrow 0}\frac{1}{\delta}(u_{\varepsilon,\delta}^h(x,t)-K_1(x))^- dxdt.
\end{equation}
\begin{equation}\label{3.5}
\xi^h(dx,dt)=\lim_{\delta\rightarrow 0}\lim_{\varepsilon\rightarrow 0}\frac{1}{\varepsilon}(u_{\varepsilon,\delta}^h(x,t)-K_2(x))^+ dxdt.
\end{equation}
Now we show that the measures $\eta^h,\xi^h$ are absolutely continuous w.r.t. the Lebesgue measure $dx\times dt$. Observe that the following equation holds:
\begin{eqnarray}\label{3.6}
       &&\frac{\partial{(u_{\varepsilon,\delta}^h(x,t)-K_1(x))}}{\partial{t}}\nonumber\\
       &=&\frac{\partial^2{(u_{\varepsilon,\delta}^h(x,t)-K_1(x))}}{\partial{x^2}}
       +f(u_{\varepsilon,\delta}^h(x,t))
       +\sigma(u_{\varepsilon,\delta}^h(x,t))\dot{h}(x,t)\nonumber\\
       &&-\alpha (u_{\varepsilon,\delta}^h(x,t)-K_1(x))\nonumber\\
       && +\frac{1}{\delta}(u_{\varepsilon,\delta}^h(x,t)-K_1(x))^- -\frac{1}{\varepsilon}(u_{\varepsilon,\delta}^h(x,t)-K_2(x))^+ \\
       && +\frac{\partial^2{K_1(x)}}{\partial{x^2}}-\alpha K_1(x).
      \nonumber
\end{eqnarray}
Multiplying both sides of the above equation by $(u_{\varepsilon,\delta}^h(x,t)-K_1(x))^-$ and integrating against $dx$ we obtain
\begin{eqnarray}\label{3.7}
       &&-\frac{\partial}{\partial t}\int_0^1[(u_{\varepsilon,\delta}^h(x,t)-K_1(x))^{-}]^2dx\nonumber\\
       &=&\int_0^1|\frac{\partial{(u_{\varepsilon,\delta}^h(x,t)-K_1(x))^-}}{\partial{x}}|^2dx
       +\alpha\int_0^1[(u_{\varepsilon,\delta}^h(x,t)-K_1(x))^-]^2dx\nonumber\\
       &&+\int_0^1 f(u_{\varepsilon,\delta}^h(x,t))(u_{\varepsilon,\delta}^h(x,t)-K_1(x))^-dx\nonumber\\
       &&
       +\int_0^1 \sigma(u_{\varepsilon,\delta}^h(x,t))(u_{\varepsilon,\delta}^h(x,t)-K_1(x))^- \dot{h}(x,t)dx\nonumber\\
       && +\frac{1}{\delta}\int_0^1[(u_{\varepsilon,\delta}^h(x,t)-K_1(x))^-]^2dx -\frac{1}{\varepsilon}\int_0^1(u_{\varepsilon,\delta}^h(x,t)-K_2(x))^+(u_{\varepsilon,\delta}^h(x,t)-K_1(x))^-dx \nonumber\\
       &&+\int_0^1\frac{\partial^2{K_1(x)}}{\partial{x^2}}(u_{\varepsilon,\delta}^h(x,t)-K_1(x))^-dx
       -\alpha \int_0^1K_1(x)(u_{\varepsilon,\delta}^h(x,t)-K_1(x))^-dx.
\end{eqnarray}
Applying the chain rule and integrating w.r.t. $t$ from $0$ to $T$ yield
\begin{eqnarray}\label{3.8}
       &&\int_0^1[(u_{\varepsilon,\delta}^h(x,T)-K_1(x))^{-}]^2dx
       +\frac{1}{\delta}\int_0^Te^{-\alpha (T-t)}\int_0^1[(u_{\varepsilon,\delta}^h(x,t)-K_1(x))^-]^2dxdt \nonumber\\
       &&+\int_0^Te^{-\alpha (T-t)}\int_0^1|\frac{\partial{(u_{\varepsilon,\delta}^h(x,t)-K_1(x))^-}}{\partial{x}}|^2dxdt\nonumber\\
       &=& -\int_0^Te^{-\alpha (T-t)}\int_0^1 f(u_{\varepsilon,\delta}^h(x,t))(u_{\varepsilon,\delta}^h(x,t)-K_1(x))^-dxdt\nonumber\\
       && -\int_0^Te^{-\alpha (T-t)}\int_0^1\sigma(u_{\varepsilon,\delta}^h(x,t))\dot{h}(x,t)(u_{\varepsilon,\delta}^h(x,t)-K_1(x))^-dxdt\nonumber\\
       &&+\frac{1}{\varepsilon}\int_0^Te^{-\alpha (T-t)}\int_0^1(u_{\varepsilon,\delta}^h(x,t)-K_2(x))^+
       (u_{\varepsilon,\delta}^h(x,t)-K_1(x))^-dxdt \nonumber\\
       &&-\int_0^Te^{-\alpha (T-t)}\int_0^1\frac{\partial^2{K_1(x)}}{\partial{x^2}}(u_{\varepsilon,\delta}^h(x,t)-K_1(x))^-dxdt\nonumber\\
       &&+\alpha \int_0^Te^{-\alpha (T-t)}\int_0^1 K_1(x)(u_{\varepsilon,\delta}^h(x,t)-K_1(x))^-dxdt.
\end{eqnarray}
Note that
$$(u_{\varepsilon,\delta}^h(x,t)-K_2(x))^+
       (u_{\varepsilon,\delta}^h(x,t)-K_1(x))^-\leq 0,$$
and for any $\varepsilon>0$, there exists a constant $C_{\varepsilon}$ such that
\begin{equation}\label{3.8-1}
ab\leq \varepsilon \frac{1}{\delta} a^2+ C_{\varepsilon} \delta b^2, \quad \quad a, b\in R
\end{equation}
From (\ref{3.8}) and (\ref{3.8-1}) we deduce that
 \begin{eqnarray}\label{3.9}
       &&\frac{1}{\delta}\int_0^Te^{-\alpha (T-t)}\int_0^1[(u_{\varepsilon,\delta}^h(x,t)-K_1(x))^-]^2dxdt \nonumber\\
       &\leq &\frac{1}{2}\frac{1}{\delta}\int_0^Te^{-\alpha (T-t)}\int_0^1[(u_{\varepsilon,\delta}^h(x,t)-K_1(x))^-]^2dxdt\nonumber\\
       &&+C\delta \int_0^Te^{-\alpha (T-t)}\int_0^1\sigma(u_{\varepsilon,\delta}^h(x,t))^2\dot{h}^2(x,t)dxdt\nonumber\\
       && +C\delta \int_0^Te^{-\alpha (T-t)}\int_0^1 f(u_{\varepsilon,\delta}^h(x,t))^2dxdt +C \delta \int_0^Te^{-\alpha (T-t)}\int_0^1(\frac{\partial^2{K_1(x)}}{\partial{x^2}})^2dxdt\nonumber\\
       &&+C \delta \int_0^Te^{-\alpha (T-t)}\int_0^1(K_1(x))^2dxdt.
\end{eqnarray}
In particular, we obtain that
$$\frac{1}{\delta^2}\int_0^Te^{-\alpha (T-t)}\int_0^1[(u_{\varepsilon,\delta}^h(x,t)-K_1(x))^-]^2dxdt\leq C_{T,h}, $$
Since $f, \sigma$ are bounded, we see that  if $||h||_2^{\infty}=\int_0^{\infty}\int_0^1 \dot{h}^2(x,s)dxds<\infty $, then
$$C_{T, h}\leq C(1+||h||_2^{\infty})$$
for some constant $C$ independent of $T$.
Subtracting a weak convergent subsequence if necessary, the above inequality  together with (\ref{3.4}) implies that $\eta^h$ is absolutely continuous w.r.t. $dxdt$ and
\begin{eqnarray}\label{3.10}
 &&\int_0^Te^{-\alpha (T-t)}\int_0^1(\dot{\eta}^h(x,t))^2dxdt \nonumber\\
 &\leq& \lim_{\varepsilon\rightarrow 0}\lim_{\delta\rightarrow 0}\frac{1}{\delta^2}\int_0^Te^{-\alpha (T-t)}\int_0^1 |u_{\varepsilon,\delta}^h(x,t)-K_1(x))^-|^2 dxdt\nonumber\\
 &\leq&  C_{T,h}.
 \end{eqnarray}
 The proof of the corresponding conclusion for $\xi^h$ is similar.\quad\quad $\blacksquare$
\vskip 0.2cm
Let $C^{\gamma}([0,1])$ denote the Banach space of $\gamma$-H\"{o}lder continuous functions on $[0,1]$ equipped with the
H\"{o}lder norm $||\cdot||_{\gamma}$.
 \begin{prop}
 Assume (F1). Let $u^h(x,t)$ be the solution of equation (\ref{3.1}). For $0<\gamma <\frac{1}{2}$,  we have
 \begin{equation}\label{3.10}
 ||u^h(\cdot, t )||_{\gamma} \leq C(1+\frac{1}{\sqrt{t}})(1+||h||_2^{\infty}).
 \end{equation}
 \end{prop}
\noindent\textbf{Proof.}\quad
  Since $\eta^h$, $\xi^h$ are absolutely continuous w.r.t. $dxdt$, it follows  that $u^h$ has the  following mild form:
   \begin{eqnarray}\label{3.11}
u^h(x,t)&=&\int_0^1G_t(x,y)u_0(y)dy+ \int_0^t \int_0^1G_{t-s}(x,y) f(y,u^h(y,s))dyds\nonumber\\
&+&\int_0^t \int_0^1G_{t-s}(x,y)\sigma(y,u^h(y,s))\dot{h}(y,s)dyds+\int_0^t \int_0^1G_{t-s}(x,y)\dot{\eta}^h(y,s)dyds\nonumber\\
&-&\int_0^t \int_0^1G_{t-s}(x,y)\dot{\xi}^h(y,s)dyds.
\end{eqnarray}
Where $G_t(x,y)=e^{-\alpha t}P_t(x,y)$ and $P_t(x,y)$ is the heat kernel of the Neumann Laplacian on $[0,1]$.
The Proposition will follow  if we prove that each of the five terms on the right has the bound (\ref{3.10}). Recall the following inequality proved in \cite{W}: for $0<\gamma<\frac{1}{2}$,
\begin{equation}\label{3.11-1}
\left (\int_0^{\infty}\int_{0}^1|P_u(x_1,y)-P_u(x_2,y)|^2dudy\right )^{\frac{1}{2}}\leq C |x_1-x_2|^{\gamma}.
\end{equation}
 By the property of the heat kernel, it holds that
$$ |\int_0^1G_t(x_1,y)u_0(y)dy-\int_0^1G_t(x_2,y)u_0(y)dy|\leq Ce^{-\alpha t}\frac{1}{\sqrt{t}} |x_1-x_2|.$$
The remaining  terms can be treated in a similar way.   Let us  look at one of the terms, say, the fourth term $F(x,h):=\int_0^t \int_0^1G_{t-s}(x,y)\dot{\eta}^h(y,s)dyds$. By Theorem 4.1, we have
\begin{eqnarray}\label{3.12}
|F(x,t)|& \leq & (\int_0^t e^{-\alpha(t-s)}\int_0^1P_{t-s}^2(x,y)dyds)^{\frac{1}{2}}(\int_0^t e^{-\alpha(t-s)} \int_0^1(\dot{\eta}^h)^2(y,s)dyds)^{\frac{1}{2}}\nonumber\\
&\leq& C (1+ ||h||_2^{\infty}).
\end{eqnarray}
 For $0<\gamma< \frac{1}{2}$, $x_1, x_2 \in [0,1]$, we have
\begin{eqnarray}\label{3.13}
&&|F(x_1,t)-F(x_2,t)|\nonumber\\
&\leq & (\int_0^t e^{-\alpha(t-s)}\int_0^1(P_{t-s}(x_1,y)- P_{t-s}(x_2,y))^2dyds)^{\frac{1}{2}}(\int_0^t e^{-\alpha(t-s)} \int_0^1(\dot{\eta}^h)^2(y,s)dyds)^{\frac{1}{2}}\nonumber\\
&\leq& C[1+(||h||_2^{\infty})^{\frac{1}{2}}]   |x_1-x_2|^{\gamma},
\end{eqnarray}
where $(\ref{3.11-1})$ and (\ref{3.2-1}) have been used. Combining (\ref{3.12}) and (\ref{3.13}) we conclude
$$||F(\cdot, t)||_{\gamma}\leq C(1+||h||_2^{\infty}).$$
The proof is complete.\quad\quad $\blacksquare$
\vskip 0.4cm

Let $v(\cdot, \cdot )\in C([0,1]\times R)$. For $ t_1<t_2$, define
\begin{equation}\label{3.14}
I_{t_1}^{t_2}(v)=\inf\{\frac{1}{2}|\dot{h}|^2_{L^2([0,1]\times [t_1,t_2])}; v=u^h\},
\end{equation}
where $u^h$ is the solution of the equation (\ref{3.1}) on the time interval $[t_1, t_2]$.
Introduce
\begin{equation}\label{3.15}
E=\{z\in C([0,1]); \quad K_1(x)\leq z(x)\leq K_2(x), x\in [0,1]\}.
\end{equation}
$E$ is a complete metric space  equipped with the metric deduced from the  maximum norm $||\cdot ||_{\infty}$
on $C([0,1])$.
Let $s>0$, $t>0$. Set
$$K_{0,t}(s)=\{ \phi\in C([0,t]; E); I_0^t(\phi)\leq s\},$$
and
$$K_{z,0,t}(s)=\{ \phi\in C([0,t]; E); \phi(0)=z, I_0^t(\phi)\leq s\}.$$
For $z\in E$, define
\begin{equation}\label{3.16}
J(z)= \inf\{I_{0}^{t}(v); v\in C([0,t]; E), v(\cdot, 0)=0, v(\cdot, t)=z, t>0 \}.
\end{equation}
\begin{theorem}
We have
\begin{equation}\label{3.17}
J(z)= \inf\{I_{-\infty }^{0}(v); v(\cdot, 0)=z, \lim_{t\rightarrow -\infty}v(\cdot, t)=0 \}
\end{equation}
\end{theorem}
\vskip 0.3cm
\noindent{\bf Proof}. Let $t>0$ and $v\in C([0,t]; E)$ with $v(0)=0$, $v(t)=z$. Define
$$\bar{v}(s)=\left \{\begin{array}{ll} v(s+t) & if \quad s\in [-t, 0],\\
0& if \quad s\leq -t.\end{array}\right. $$
Then $\bar{v}(0)=z$, $\lim_{s\rightarrow -\infty}\bar{v}(\cdot, s)=0$ and $I_{-\infty }^{0}(\bar{v})=I_{0}^{t}(v)$. Consequently,
$$ \inf\{I_{-\infty }^{0}(v); v(\cdot, 0)=z, \lim_{t\rightarrow -\infty}v(\cdot, t)=0 \}\leq I_{-\infty }^{0}(\bar{v})=I_{0}^{t}(v)$$
As $t$, $v$ are arbitrary, we deduce that
$$J(z)\geq \inf\{I_{-\infty }^{0}(v); v(\cdot, 0)=z, \lim_{t\rightarrow -\infty}v(\cdot, t)=0 \}.$$
To prove the opposite inequality, we may assume
$$\inf\{I_{-\infty }^{0}(v); v(\cdot, 0)=z, \lim_{t\rightarrow -\infty}v(\cdot, t)=0 \}<\infty .$$
In this case, following the same method as in \cite{S}, \cite{CR}) we can show that the $\inf$ can be attained, i.e.,  there exists $v_0$ with $v_0(\cdot, 0)=z$, $\lim_{t\rightarrow -\infty}v_0(\cdot, t)=0$ such that
\begin{equation}\label{3.18}
I_{-\infty}^0(v_0)= \inf\{I_{-\infty }^{0}(v); v(\cdot, 0)=z, \lim_{t\rightarrow -\infty}v(\cdot, t)=0 \}.
\end{equation}
In view of the assumptions on $K_1(x)$ and $K_2(x)$, there exists $\varepsilon_0>0$ such that
\begin{equation}\label{3.19}
K_1(x)< -\varepsilon_0 <0<\varepsilon_0<K_2(x).
\end{equation}
For any $\varepsilon>0$, there exists $n_0$ such that  $||v_0(\cdot, t)||_{\infty}\leq \varepsilon_0\wedge \varepsilon$ for $t\leq -n_0+2 $ and $I_{-n_0}^{-n_0+1}(v_0)\leq \varepsilon$. Write $v_0(t)$ for $v_0(\cdot,t)$ and define
$$v_{n_0}(t)=\left \{\begin{array}{ll} v_0(t) & if \quad t\in [-n_0+1, 0],\\
(t+n_0)v_0(t)& if \quad -n_0\leq t\leq -n_0+1.\end{array}\right. $$
Set $\bar{v}_{n_0}(t)={v}_{n_0}(t-n_0)$, for $0\leq t\leq n_0$. Then $\bar{v}_{n_0}(0)=0$, $\bar{v}_{n_0}(n_0)=z$.
For $v(\cdot, \cdot )\in C([0,1]\times R)$ and $ t_1<t_2$, define
\begin{equation}\label{3.20}
S_{t_1}^{t_2}(v)=\inf\{\frac{1}{2}|\dot{h}|^2_{L^2([0,1]\times [t_1,t_2])}; v=v^h\},
\end{equation}
where $v^h$ is the solution of the following skeleton equation:
\begin{eqnarray}\label{3.21}\left\{
\begin{array}{l}
\frac{\partial v^h(x,t)}{\partial t}-\frac{\partial^2 v^h(x,t)}{\partial x^2}+\alpha v^h(x,t)=f(x,v^h)+\sigma(x,v^h)\dot{h}(x,t);\\
u^h(\cdot,0)=u_0;\\
\frac{\partial u^h}{\partial x}(0,t)=\frac{\partial u^h}{\partial x}(1,t)=0.
\end{array}\right.
\end{eqnarray}
If $||v(t)||_{\infty}\leq \varepsilon_0$ for $t_1\leq t\leq t_2$, it is clear that $ S_{t_1}^{t_2}(v)=I_{t_1}^{t_2}(v)$ because in this case, the extra forces $\eta^h$, $\xi^h$ in (\ref{3.1}) disappear. From the proof of Proposition 7
in \cite{S}, we know that there exists a constant $C$ such that
$$S_{-n_0}^{-n_0+1}(v_{n_0})\leq C [S_{-n_0}^{{-n_0+1}}(v_{0})+\sup_{-n_0\leq t\leq -n_0+1}(||v_0(t)||_{\infty})^2].$$
Taking into account of the choice of $n_0$ it follows that
\begin{eqnarray}\label{3.22}
&&J(z)\leq I_0^{n_0}(\bar{v}_{n_0})=I_{-n_0}^0(v_{n_0})\nonumber\\
&\leq& I_{-n_0+1}^0(v_{n_0})+I_{-n_0}^{-n_0+1}(v_{n_0})\nonumber\\
&\leq& I_{-\infty}^0(v_{0})+S_{-n_0}^{{-n_0+1}}(v_{n_0})\nonumber\\
&\leq& I_{-\infty}^0(v_{0})+C [S_{-n_0}^{{-n_0+1}}(v_{0})+\sup_{-n_0\leq t\leq -n_0+1}(||v_0(t)||_{\infty})^2]\nonumber\\
&\leq &I_{-\infty}^0(v_{0})+C [I_{-n_0}^{{-n_0+1}}(v_{0})+\sup_{-n_0\leq t\leq -n_0+1}(||v_0(t)||_{\infty})^2]\nonumber\\
&\leq& I_{-\infty}^0(v_{0})+C (\varepsilon+\varepsilon^2),
\end{eqnarray}
where $C$ is an independent constant.
Since $\varepsilon$ is arbitrary, this proves
$$J(z)\leq  \inf\{I_{-\infty }^{0}(v); v(\cdot, 0)=z, \lim_{t\rightarrow -\infty}v(\cdot, t)=0 \}$$
by the choice of $v_0$.\quad\quad $\blacksquare$
\begin{prop}
The functional $J(\cdot): E\rightarrow [0, +\infty]$ is lower semi-continuous with compact level sets, i.e.,
for $r\geq 0$, $K(r)=\{z\in E; J(z)\leq r\}$ is compact.
\end{prop}
The proof of this proposition is very  similar to that of Theorem 5.5 in \cite{CR} (see also Section 6 in \cite{S}). We omit the details.

\section{Exponential tightness of invariant measures}
\setcounter{equation}{0}

Let $\mu_{\varepsilon}$ denote the unique invariant probability measure of the solution $u_{\varepsilon}(t), t\geq 0$ of equation (\ref{I1}).
The following result is an exponential tightness for the invariant measures.

\begin{theorem}
Assume (F1). For any $L>0$, there exists a compact subset $K_L\subset C([0,1])$ such that
\begin{equation}\label{ex1}
\mu_{\varepsilon}(K_L^c)\leq exp(-\frac{L}{\varepsilon^2})
\end{equation}
for $\varepsilon\leq \varepsilon_0$, where $\varepsilon_0$ is a positive constant.
\end{theorem}

\noindent\textbf{Proof.}\quad By the invariance of $\mu_{\varepsilon}$, we have
\begin{equation}
\mu_{\varepsilon}(K_L^c)=\int_EP(u_{\varepsilon}(\cdot, 1, g)\in K_L^c)\mu_{\varepsilon}(dg).
\end{equation}
Thus it is sufficient to prove that there exists a compact subset $K_L\subset C([0,1])$ such that
\begin{equation}
P(u_{\varepsilon}(\cdot, 1, g)\in K_L^c)\leq exp(-\frac{L}{\varepsilon^2})
\end{equation}
for $\varepsilon\leq \varepsilon_0$ and  all $ g\in{C}([0,1])$ with $ K_1\leq{}g\leq{}K_2$.
Put
\begin{eqnarray}\label{2.2}
        v_{\varepsilon}(x,t,g)&=&\int_0^t\int_{0}^1G_{t-s}(x,y)f(y, u_{\varepsilon}(y,s,g))dyds \nonumber\\
        &&+\varepsilon\int_0^t\int_{0}^1G_{t-s}(x,y)\sigma(y, u_{\varepsilon}(y,s,g))W(dy,ds).
\end{eqnarray}
 Then $u_{\varepsilon}$ can be written in the form
\begin{eqnarray*}
        u_{\varepsilon}(x,t,g)-\int_{0}^1G_t(x,y)g(y)dy&=&v_{\varepsilon}(x,t,g)+\int_0^t\int_{0}^1G_{t-s}(x,y)\eta_{\varepsilon}(g)(dx,dt)\\
        &&-\int_0^t\int_{0}^1G_{t-s}(x,y)\xi_{\varepsilon}(g)(dx,dt),
\end{eqnarray*}
where $\eta_{\varepsilon}(g)$, $\xi_{\varepsilon}(g)$ indicates the dependence of the random measures on the initial condition $g$.
Put
$$\bar{u}_{\varepsilon}(x,t,g)=u_{\varepsilon}(x,t,g)-\int_{0}^1G_t(x,y)g(y)dy.$$
Then $(\bar{u}_{\varepsilon}, \eta_{\varepsilon}, \xi_{\varepsilon} )$ solves a random obstacle problem (\ref{1.1}) with $v(x,t)$ replaced by $v_{\varepsilon}(x,t,g)$.
 As shown in \cite{YZ}, there exists a  continuous functional $\Phi_1$  from ${C}([0,1]\times[0,1])$ to ${C}([0,1])$ such that  $\bar{u}_{\varepsilon}(\cdot,1,g)=\Phi_1(v_{\varepsilon}(\cdot,g))$, where $v_{\varepsilon}(\cdot,g)=v_{\varepsilon}(\cdot,\cdot,g)$.
 Since $f(u_{\varepsilon}(y,s,g))$, $\sigma(u_{\varepsilon}(y,s,g))$ are bounded by a constant independent of $g$, by Proposition 4 in \cite{S}, for $L>0$ there exists a compact subset $K_L^{\prime}\subset {C}([0,1]\times[0,1])$ such that
\begin{equation}\label{2.3}
P(v_{\varepsilon}(\cdot, \cdot, g)\in (K_L^{\prime})^c)\leq exp(-\frac{L}{\varepsilon^2})
\end{equation}
for $\varepsilon\leq \varepsilon_0$ and  all $ g\in{C}([0,1])$ with $ K_1\leq{}g\leq{}K_2$.
 On the other hand, it is easy to see that there is a compact subset $K_0\subset C([0,1])$ such that
$$ \{ \int_{0}^1G_1(x,y)g(y)dy; \quad K_1\leq{}g\leq{}K_2 \}\subset K_0.$$
Let $\Phi_1(K_L^{\prime})$ denote the image of $K_L^{\prime}$ under the map $\Phi_1$. Put $K_L=\Phi_1(K_L^{\prime})+K_0$. Then we have
\begin{equation}\label{2.4}
P(u_{\varepsilon}(\cdot, 1, g)\in K_L^c)\leq P(v_{\varepsilon}(\cdot, \cdot, g)\in (K_L^{\prime})^c)\leq exp(-\frac{L}{\varepsilon^2})
\end{equation}
for $\varepsilon\leq \varepsilon_0$ and  all $ g\in{C}([0,1])$ with $ K_1\leq{}g\leq{}K_2$.
 This finishes the proof.\quad\quad $\blacksquare$ \\
\section{Statement of large deviations}
\setcounter{equation}{0}
The following result is a large deviation principle of $u_{\varepsilon}$ ( the solution of the equation (\ref{I1})) on the path space $C([0,1]\times[0,T])$. The proof of the theorem is very similar to that of Theorem 5.1 in \cite{XZ} where a large deviation principle was proved for SPDEs with reflection at $0$. We just need to use Theorem 7 in \cite{BDM} to improve  Theorem 5.1 in \cite{XZ} to a uniform large deviation principle on compact sets.

\begin{theorem} Assume (F1). Then, the laws $\nu_{\varepsilon}^g$ of $\{u_\varepsilon(\cdot,
\cdot, g)\}_{\varepsilon>0}$ satisfy a large deviation principle on
$C([0,1]\times[0,T])$ with the rate function $I_0^T(\cdot)$ uniformly on compact sets , i.e., given any compact subset $K$. we have

(i) for any closed subset $C\subset C([0,1]\times[0,T])$,

$$\limsup_{\varepsilon \rightarrow 0}\varepsilon^2 \log \sup_{g\in K}\nu_{\varepsilon}^g(C)\leq -\inf_{f\in C}I_0^T(f).$$

(ii) for any open subset  $G\subset C([0,1]\times[0,T]) $,

$$\liminf_{\varepsilon \rightarrow 0}\varepsilon^2 \log \inf_{g\in K}\nu^g_{\varepsilon}(G)\geq -\inf_{f\in
G}I_0^T(f).$$
\end{theorem}
\vskip 0.4cm

Let $\mu_{\varepsilon}$ denote the unique invariant probability measure of the solution $u_{\varepsilon}(t), t\geq 0$.
Introduce the following assumption:

\noindent{\bf (H)}. Assume $f(x,0)=0$ and that there exists a constant  $c<\alpha$ such that
\begin{equation}\label{6.0}
|f(x,u_1)-f(x,u_2)|\leq c |u_1-u_2|, u_1, u_2\in R.
\end{equation}
\vskip 0.3cm

Here is the main result of this paper:

\begin{theorem} Suppose the conditions (F1), (F2) and (H) hold. Then $\mu_{\varepsilon}, \varepsilon>0$ satisfies a large deviation principle on
$E$ with the rate function $J(\cdot)$, i.e.,

(i) for any closed subset $C\subset E$,

$$\limsup_{\varepsilon \rightarrow 0}\varepsilon^2 \log \mu_{\varepsilon}(C)\leq -\inf_{z\in C}J(z).$$

(ii) for any open subset  $G\subset E$,

$$\liminf_{\varepsilon \rightarrow 0}\varepsilon^2 \log \mu_{\varepsilon}(G)\geq -\inf_{z\in
G}J(z).$$
\end{theorem}
To prove this theorem, it is well known (see e.g. \cite{F}, \cite{S} and \cite{CR}) that it suffices to establish the following:

1. lower bounds: for any $\delta>0, \gamma>0$ and $z^*\in E$ there exists $\varepsilon_0>0$ such that
\begin{equation}\label{6.1}
\mu_{\varepsilon}(\{z\in E: ||z-z^*||_{\infty}<\delta \})\geq exp(-\frac{J(z^*)+\gamma}{\varepsilon^2}), \quad \varepsilon \leq \varepsilon_0,
\end{equation}

2. upper bounds: for any $s\geq 0$, $\delta>0$, $\gamma>0$ there exists $\varepsilon_0>0$ such that
\begin{equation}\label{6.2}
\mu_{\varepsilon}(\{z\in E; dis_E(z, K(s))\geq \delta\})\leq exp(-\frac{s-\gamma}{\varepsilon^2}), \quad \varepsilon\leq \varepsilon_0,
\end{equation}
where $K(s)=\{ z\in E; J(z)\leq s\}.$

These will be proved in  Section 7 and 8.

\section{Lower bounds}
\setcounter{equation}{0}
Consider
\begin{eqnarray}\label{4.1}\left\{
\begin{array}{l}
\frac{\partial u^0(x,t)}{\partial t}=\frac{\partial^2 u^0(x,t)}{\partial x^2}-\alpha u^0(x,t)+f(x,u^0)+\eta^z-\xi^z;\\
K_1(x)\leq u^0(x,t)\leq K_2(x);\\
\int_0^T\int_0^1(u^0(x,t)-K_1(x))\eta^z(dx,dt)=\int_0^T\int_0^1(K_2(x)-u^0(x,t))\xi^z(dx,dt)=0\\
u^0(\cdot,0)=z;\\ \frac{\partial u^0}{\partial x}(0,t)=\frac{\partial u^0}{\partial x}(1,t)=0.
\end{array}\right.
\end{eqnarray}
Write the solution of (\ref{4.1}) as $u^0(z,x,t)$ to stress the dependence on the initial function $z$. Denote by $B$ the Banach space $C([0,1])$ with the maximum norm $||\cdot ||\_{\infty}$.
\begin{lemma}
Assume the conditions (F1), (F2) and (H).
Then it holds that
\begin{equation}\label{4.3}
||u^0(z,\cdot ,t)||_{\infty}\leq e^{-\alpha_1 t}||z||_{\infty},
\end{equation}
where $\alpha_1=\alpha-c>0$.
\end{lemma}
\noindent\textbf{Proof.}\quad
Set
$$A=\frac{\partial^2}{\partial x^2}-\alpha I.$$
We write $f(g)(x)$ for $f(x,g(x))$ for brevity. A similar notation will be used for $\sigma(x,g(x))$.
First we claim that for any $g\in E\cap D(A)$, there exists
$$l_g\in \partial ||g||_B:=\{l\in B^*; ||l||_{B^*}=1, <l,g>=||g||_B\}$$
such that
\begin{equation}\label{4.4}
<Ag, l_g>+<f(g),l_g>\leq -\alpha_1 ||g||_B.
\end{equation}
In fact, choose $l_g=\delta_{x_0}$ if $g(x_0)=max_{x\in [0,1]}|g(x)|$  and
  $l_g=-\delta_{x_0}$ if $g(x_0)=-max_{x\in [0,1]}|g(x)|$. Then
 \begin{eqnarray}\label{4.5}
 <Ag, l_g>+<f(g),l_g>&\leq &-\alpha ||g||_B +f(x_0, g(x_0))\nonumber\\
 &\leq& (c-\alpha )||g||_B.
 \end{eqnarray}
 Furthermore, if $g(x_0)=max_{x\in [0,1]}|g(x)|$, then, as $g\in E$,
 $$<(g-K_1)^-,l_g>=(g(x_0)-K_1(x_0))^-=0.$$
 If $g(x_0)=-max_{x\in [0,1]}|g(x)|$, then
 $$ <(g-K_1)^-,l_g>=-(g(x_0)-K_1(x_0))^-\leq 0,$$
 $$ <(g-K_2)^+,l_g>=0.$$
 So in all cases we have
 \begin{eqnarray}\label{4.6}
 &&<Ag, l_g>+<f(g),l_g>\nonumber\\
 &&+\frac{1}{\delta}<(g-K_1)^-,l_g>-\frac{1}{\varepsilon}<(g-K_2)^+,l_g>\nonumber\\
 &\leq& (c-\alpha )||g||_B,
 \end{eqnarray}
 for all $\varepsilon >0$, $\delta>0$.
 Consider the approximating equations:
\begin{eqnarray}\label{4.7}
       \frac{\partial{u_{\varepsilon,\delta}^0(x,t)}}{\partial{t}}&=&
       \frac{\partial^2{u_{\varepsilon,\delta}^0(x,t)}}{\partial{x^2}}
       +f(u_{\varepsilon,\delta}^0(x,t))-\alpha u_{\varepsilon,\delta}^0(x,t)\nonumber\\
       && +\frac{1}{\delta}(u_{\varepsilon,\delta}^0(x,t)-K_1(x))^- -\frac{1}{\varepsilon}(u_{\varepsilon,\delta}^0(x,t)-K_2(x))^+; \\
       u_{\varepsilon,\delta}^0(x,0)&=&u_0(x); \nonumber \\
       \frac{\partial u_{\varepsilon,\delta}^0}{\partial x}(0,t)&=&\frac{\partial u_{\varepsilon,\delta}^0}{\partial x}(1,t)=0.\nonumber
\end{eqnarray}
By the chain rule, we have
\begin{eqnarray}\label{4.8}
       \frac{d^-}{dt}||u_{\varepsilon,\delta}^0(t)||_B&\leq &
       <Au_{\varepsilon,\delta}^0(t), l_{u_{\varepsilon,\delta}^0(t)}>
       +<f(u_{\varepsilon,\delta}^0(t)),l_{u_{\varepsilon,\delta}^0(t)}> \nonumber\\
       && +\frac{1}{\delta}<u_{\varepsilon,\delta}^0(t)-K_1)^-, l_{u_{\varepsilon,\delta}^0(t)}> -\frac{1}{\varepsilon}<u_{\varepsilon,\delta}^0(t)-K_2)^+, l_{u_{\varepsilon,\delta}^0(t)}> \nonumber\\
       &\leq & -\alpha_1 ||u_{\varepsilon,\delta}^0(t)||_B.
\end{eqnarray}
This yields
$$||u_{\varepsilon,\delta}^0(t)||_B\leq e^{-\alpha_1t}||z||_B,$$
where $z$ is the initial function. Because the constants on the right side are independent of $\delta, \varepsilon$, let $\delta\rightarrow 0$ and $\varepsilon\rightarrow 0$ to get (\ref{4.3}).\quad\quad $\blacksquare$
\vskip 0.4cm
Fix $z^{\ast}\in E$ with $J(z^{\ast})<\infty$. For any $\gamma >0$, by the definition of $J(z^*)$ there exists a function $\psi$ and $T_0>0$ such that $\psi (0)=0$, $\psi (T_0)=z^*$ and $\psi=u^{\bar{h}}$ for some $\bar{h}$ with
$\frac{1}{2}|\dot{\bar{h}}|^2_{L^2([0,1]\times [0,T_0])}<J(z^{\ast})+\frac{\gamma}{2}$.
Define for $T>0$,
\begin{equation}\label{4.9}
\dot{h}^T(x,t)=\left\{\begin{array}{ll}
0, & (x,t)\in [0,1]\times [0,T],\\
\dot{\bar{h}}(x,t-T)& (x,t)\in [0,1]\times [T,T+T_0].
\end{array}\right.
\end{equation}
Consider the following PDE with reflection:
\begin{eqnarray}\label{4.10}
\frac{\partial u^T(z, x,t)}{\partial t}&=& \frac{\partial^2 u^T(z,x,t)}{\partial x^2}+f(x,u^T)+\sigma(x,u^T)\dot{h}^T(x,t)\nonumber\\
&& +\eta^T-\xi^T;\\
u^T(\cdot,0)&=&z, \quad  u^h(0,t)=u^h(1,t)=0.\nonumber
\end{eqnarray}
Clearly,
\begin{equation}\label{4.11}
u^T(z, x,t)=\left\{\begin{array}{ll}
u^0(z,x,t), & (x,t)\in [0,1]\times [0,T],\\
 u^{\bar{h}}(u^0(z,\cdot,T),x,t-T)& (x,t)\in [0,1]\times [T,T+T_0].
\end{array}\right.
\end{equation}
Set $\bar{\psi}(x,t)=u^T(z, x,t+T)$ for $t\geq 0$. Then we have
\begin{eqnarray}\label{4.12}
\frac{\partial \bar{\psi}(z, x,t)}{\partial t}&=& \frac{\partial^2 \bar{\psi}(z,x,t)}{\partial x^2}+f(x,\bar{\psi})+\sigma(x,\bar{\psi})\dot{\bar{h}}(x,t)\nonumber\\
&& +\bar{\eta}-\bar{\xi};\\
\bar{\psi}(x,0)&=&u^0(z,x,T).\nonumber
\end{eqnarray}
Recall that $\psi=u^{\bar{h}}$ satisfies the following reflected PDE:
\begin{eqnarray}\label{4.13}
\frac{\partial \psi (x,t}{\partial t}&=& \frac{\partial^2 \psi(x,t)}{\partial x^2}+f(x,\psi )+\sigma(x,\psi )\dot{\bar{h}}(x,t)\nonumber\\
&& +\eta-\xi ;\\
\psi(x,0)&=&0.\nonumber
\end{eqnarray}
Put
$$F(t):=\sup_{0\leq s\leq t}\sup_{x\in [0,1]}|\psi(x,s)-\bar{\psi}(x,s)|=\sup_{0\leq s\leq t}||\psi(s)-\bar{\psi}(s)||_{\infty}.$$
We have the following result:
\begin{prop} Assume (F1). We have
\begin{equation}\label{4.14}
F(T_0)\leq C(T_0, ||\bar{h}||)\sup_{x\in [0,1]}|u^0(z,x,T)|.
\end{equation}
\end{prop}
\noindent\textbf{Proof.}\quad
Set
\begin{eqnarray}\label{4.15}
A(x,t)&=&\int_0^t \int_0^1G_{t-s}(x,y) f(y,\psi(y,s))dyds\nonumber\\
&&+\int_0^t \int_0^1G_{t-s}(x,y)\sigma(y,\psi(y,s))\dot{\bar{h}}(y,s)dyds.
\end{eqnarray}
\begin{eqnarray}\label{4.16}
\bar{A}(x,t)&=&P_t(u^0(z,\cdot,T))(x)+\int_0^t \int_0^1G_{t-s}(x,y) f(y,\bar{\psi}(y,s))dyds\nonumber\\
&&+\int_0^t \int_0^1G_{t-s}(x,y)\sigma(y,\bar{\psi}(y,s))\dot{\bar{h}}(y,s)dyds.
\end{eqnarray}
Then it follows from Theorem 3.1  in \cite{YZ1} (also see Remark 3.1) that for $t\leq T_0$,
\begin{eqnarray}\label{4.17}
F(t)&\leq &2 \sup_{0\leq s\leq t}\sup_{x\in [0,1]}|A(x,s)-\bar{A}(x,s)|\nonumber\\
&\leq &2\sup_{0\leq s\leq t}\sup_{x\in [0,1]}|P_s(u^0(z,\cdot,T))(x)|\nonumber\\
&+&2\sup_{0\leq s\leq t}\sup_{x\in [0,1]}|\int_0^s \int_0^1G_{s-u}(x,y)[f(y,\psi(y,u))-f(y,\bar{\psi}(y,u))]dydu|\nonumber\\
&+& 2\sup_{0\leq s\leq t}\sup_{x\in [0,1]}|\int_0^s \int_0^1G_{s-u}(x,y)[\sigma(y,\psi(y,u))-\sigma(y,\bar{\psi}(y,u))]\dot{\bar{h}}(y,u) dydu|.\nonumber
\end{eqnarray}
\begin{eqnarray}\label{4.18}
&\leq &C_{T_0}\sup_{x\in [0,1]}|u^0(z,\cdot,T)(x)|\nonumber\\
&&+2C\sup_{0\leq s\leq t}\sup_{x\in [0,1]}|\int_0^s \sup_{y\in [0,1]}|\psi(y,u)-\bar{\psi}(y,u)| \int_0^1G_{s-u}(x,y)dydu|\nonumber\\
&&+2C\sup_{0\leq s\leq t}\sup_{x\in [0,1]}|\int_0^s \int_0^1G_{s-u}(x,y)\sup_{y\in [0,1]}|\psi(y,u)-\bar{\psi}(y,u)| |\dot{\bar{h}}(y,u)| dydu\nonumber\\
&\leq& C_{T_0}\sup_{x\in [0,1]}|u^0(z,\cdot,T)(x)|+C\int_0^t F(u)du \nonumber\\
&&+2C\sup_{0\leq s\leq t}\sup_{x\in [0,1]}|\int_0^s F(u)\int_0^1G_{s-u}(x,y)|\dot{\bar{h}}(y,u)|  dydu.
\end{eqnarray}
Now,
\begin{eqnarray}\label{4.19}
&&\sup_{0\leq s\leq t}\sup_{x\in [0,1]}\int_0^s F(u)\int_0^1G_{s-u}(x,y)|\dot{\bar{h}}(y,u)|dydu\nonumber\\
&\leq& \sup_{0\leq s\leq t}\sup_{x\in [0,1]}|\int_0^s F(u)du(\int_0^1G_{s-u}(x,y)^2dy)^{\frac{1}{2}}(\int_0^1|\dot{\bar{h}}(y,u)|^2dy)^{\frac{1}{2}}\nonumber\\
&\leq& C \sup_{0\leq s\leq t}\int_0^s F(u)du(\frac{1}{\sqrt{s-u}})^{\frac{1}{2}}(\int_0^1|\dot{\bar{h}}(y,u)|^2dy)^{\frac{1}{2}}\nonumber\\
&\leq& C \sup_{0\leq s\leq t}(\int_0^s F(u)^2\frac{1}{\sqrt{s-u}})du)^{\frac{1}{2}}(\int_0^{T_0}\int_0^1|\dot{\bar{h}}(y,u)|^2dydu)^{\frac{1}{2}}\nonumber\\
&\leq& C |\dot{\bar{h}}|_{L^2([0,T_0]\times [0,1])} \sup_{0\leq s\leq t}(\int_0^s F(u)^{2p}du)^{\frac{1}{2p}}( \int_0^s (\frac{1}{\sqrt{s-u}})^q du)^{\frac{1}{2q}}\nonumber\\
&\leq& C |\dot{\bar{h}}|_{L^2([0,T_0]\times [0,1])} C_{T_0,q}(\int_0^s F(u)^{2p}du)^{\frac{1}{2p}},
\end{eqnarray}
here $p>2, \frac{1}{p}+\frac{1}{q}=1$. Combining (\ref{4.17}) and (\ref{4.19}) we obtain for $t\leq T_0$,
\begin{eqnarray}\label{4.20}
F(t)^{2p}&\leq & C_{T_0, p}\sup_{x\in [0,1]}|u^0(z,\cdot,T)(x)|^{2p}+C_{T_0, p}\int_0^tF(u)^{2p}du\nonumber\\
&&+C_{T_0, p}|\dot{\bar{h}}|_{L^2([0,T_0]\times [0,1])}^{2p}\int_0^t F(u)^{2p}du
\end{eqnarray}
Applying Gronwall's inequality yields
$$F(T_0)^{2p}\leq \bar{C}(T_0, p, ||h||) \sup_{x\in [0,1]}|u^0(z,\cdot,T)(x)|^{2p},$$
giving (\ref{4.14}).\quad\quad $\blacksquare$

\vskip 0.3cm
Let $\mu_{\varepsilon}$ be the invariant measure of  the solution $u_{\varepsilon}$ of the reflected SPDE (\ref{I1}).
We have
\begin{theorem}
Let the assumptions (F1), (F2) and (H) hold. For any $\delta>0, \gamma>0$ and $z^*\in E$ there exists $\varepsilon_0>0$ such that
\begin{equation}\label{4.21}
\mu_{\varepsilon}(\{z\in E: ||z-z^*||_{\infty}<\delta \})\geq exp(-\frac{J(z^*)+\gamma}{\varepsilon^2}), \quad \varepsilon \leq \varepsilon_0.
\end{equation}
\end{theorem}
\noindent\textbf{Proof.}\quad
 Without loss of generality, we assume $J(z^{\ast})<\infty$. For $\gamma >0$,  there exists a function $\psi$ and $T_0>0$ such that $\psi (0)=0$, $\psi (T_0)=z^*$ and $\psi=u^{\bar{h}}$ for some $\bar{h}$ with
$\frac{1}{2}|\dot{\bar{h}}|^2_{L^2([0,1]\times [0,T_0])} <J(z^{\ast})+\frac{\gamma}{2}$. Combining Lemma 7.1 and Proposition 7.1, we see that there exists a sufficiently big constant $T>0$ such that
$$ ||z^*-u^T(z, T_0+T)||_{\infty}=||\psi(T_0)-u^T(z, T_0+T)||_{\infty}\leq \frac{\delta}{2}.$$
Thus for any $z\in E$ we have
\begin{eqnarray}\label{4.22}
&&||u_{\varepsilon}(z, T_0+T)-z^*||_{\infty}\nonumber\\
&\leq& ||u_{\varepsilon}(z, T_0+T)-u^T(z, T_0+T)||_{\infty}+||z^*-u^T(z, T_0+T)||_{\infty}\nonumber\\
&\leq & ||u_{\varepsilon}(z, T_0+T)-u^T(z, T_0+T)||_{\infty}+\frac{\delta}{2}.
\end{eqnarray}
On the other hand, using Theorem 5.1, we can find a compact subset $K$ such that for $\varepsilon \leq \varepsilon_0$,
\begin{equation}\label{exponential}
\mu_{\varepsilon}(K)\geq \frac{1}{2}.
\end{equation}
By the invariance of $\mu_{\varepsilon}$, we have, in view of (\ref{4.22}),
\begin{eqnarray}\label{4.23}
&&\mu_{\varepsilon}(\{z\in E: ||z-z^*||_{\infty}<\delta \})\nonumber\\
&=&\int_EP(||u_{\varepsilon}(z, T_0+T)-z^*||_{\infty}<\delta )\mu_{\varepsilon}(dz)\nonumber\\
&\geq& \int_EP( ||u_{\varepsilon}(z, T_0+T)-u^T(z, T_0+T)||_{\infty}<\frac{\delta}{2})\mu_{\varepsilon}(dz)\nonumber\\
&\geq& \int_EP( ||u_{\varepsilon}(z)-u^T(z)||_{C([0,1]\times [0, T_0+T])}<\frac{\delta}{2})\mu_{\varepsilon}(dz)\nonumber\\
&\geq& \int_KP( ||u_{\varepsilon}(z)-u^T(z)||_{C([0,1]\times [0, T_0+T])}<\frac{\delta}{2})\mu_{\varepsilon}(dz).
\end{eqnarray}
By the uniform large deviation principle satisfied by $u_{\varepsilon}$ (Theorem 6.1), there exists $\varepsilon_0>0$ such that
for $\varepsilon \leq \varepsilon_0$, $z\in K$,
\begin{eqnarray}\label{4.24}
P( ||u_{\varepsilon}(z)-u^T(z)||_{C([0,1]\times [0, T_0+T])}<\frac{\delta}{2})&\geq &
exp\left ( -\frac{|\dot{\bar{h}}|^2_{L^2([0,1]\times [0,T_0])}+\gamma}{\varepsilon^2}\right )\nonumber\\
&\geq&  exp(-\frac{J(z^*)+\gamma}{\varepsilon^2}), \quad \varepsilon \leq \varepsilon_0.
\end{eqnarray}
Putting (\ref{4.23}), (\ref{exponential}) and (\ref{4.24}) together, we obtain
$$\mu_{\varepsilon}(\{z\in E: ||z-z^*||_{\infty}<\delta \})\geq \frac{1}{2} exp(-\frac{J(z^*)+\gamma}{\varepsilon^2}), \quad \varepsilon \leq \varepsilon_0.$$
$\blacksquare$
\section{Upper bounds}
\setcounter{equation}{0}
\begin{lemma}
Assume (F1), (F2) and (H).
For any $\delta>0, s>0$, there exists $T_0>0$ such that
\begin{equation}\label{5.0}
\{ g(t); g\in K_{0,t}(s)\} \subset \{ z\in E; dist_E(z, K(s))\leq \frac{\delta}{2}\}, \quad t\geq T_0.
\end{equation}
\end{lemma}
The  same proof as that of Lemma 7.1 in \cite{CR} works here. We refer the reader to \cite{CR} for details.
\vskip 0.3cm

 After the necessary preparations,
 using arguments analogous to those employed  by Sowers in \cite{S} and by Cerrai, R\"ockner in \cite{CR} we can prove the following  upper bounds. Put $K(s)=\{z\in E; J(z)\leq s\}$.

\begin{theorem} Assume (F1), (F2) and (H).
For any $s\geq 0$, $\delta>0$, $\gamma>0$ there exists $\varepsilon_0>0$ such that
\begin{equation}\label{5.2}
\mu_{\varepsilon}(\{z\in E; dis_E(z, K(s))\geq \delta\})\leq exp(-\frac{s-\gamma}{\varepsilon^2}), \quad \varepsilon\leq \varepsilon_0.
\end{equation}
\end{theorem}
\noindent\textbf{Proof.}\quad
Let $L> s-\gamma$. By Theorem 5.1 there exists a compact subset $K_L\subset E$ and $\varepsilon_1>0$ such that
for $\varepsilon\leq \varepsilon_1$,
 \begin{equation}\label{5-1}
\mu_{\varepsilon}(K_L^c)\leq exp(-\frac{L}{\varepsilon^2}).
\end{equation}
By the invariance of the measure $\mu_{\varepsilon}$, for any $t\geq 0$, we have
\begin{equation}\label{5.3}
\mu_{\varepsilon}(\{z\in E; d(z, K(s))\geq \delta\})=\int_EP(d(u_{\varepsilon}(z,\cdot, t), K(s))\geq \delta )\mu_{\varepsilon}(dz).
\end{equation}
 Thus,
 \begin{eqnarray}\label{5.4}
 &&\mu_{\varepsilon}(\{z\in E; d(z, K(s))\geq \delta\})\nonumber\\
 &\leq & \mu_{\varepsilon}(K_L^c)+\int_{K_L}P(d(u_{\varepsilon}(z,\cdot, t), K(s))\geq \delta )\mu_{\varepsilon}(dz).
 \end{eqnarray}
 By Lemma 8.1, there exists $T_1>0$ such that for $t\geq T_1$,
 \begin{eqnarray}\label{5.6}
 &&P(d(u_{\varepsilon}(z,\cdot, t), K(s))\geq \delta)\nonumber\\
 &\leq& P(d(u_{\varepsilon}(z,\cdot), K_{0,t}(s))\geq \frac{\delta}{2})\nonumber\\
 &\leq & P(d(u_{\varepsilon}(z,\cdot), K_{z,0,t}(s))\geq \frac{\delta}{2}).
 \end{eqnarray}
 The uniform large deviation principle of $u_{\varepsilon}(t)$ on the path space implies that there exists $\varepsilon(t)>0$ such that
 for $z\in K_L$,
 \begin{equation}\label{5.7}
 P(d(u_{\varepsilon}(z,\cdot ), K_{z,0,t}(s))\geq \frac{\delta}{2})\leq exp(-\frac{s-\gamma/2}{\varepsilon^2}), \quad \varepsilon\leq \varepsilon(t).
 \end{equation}
 Choosing $t=T_1$ we obtain
 \begin{equation}\label{5.7}
 \int_{K_L}P(d(u_{\varepsilon}(z,\cdot, t), K(s))\geq \delta )\mu_{\varepsilon}(dz)\leq exp(-\frac{s-\gamma/2}{\varepsilon^2}),
 \end{equation}
 for $\varepsilon\leq \varepsilon_0$, where $\varepsilon_0>0$.

 \vskip 0.2cm
 Combining (\ref{5.3}), (\ref{5.4}) with (\ref{5.7}) it follows that
 $$ \mu_{\varepsilon}(\{z\in E; d(x, K(s))\geq \delta\})\leq 2 exp(-\frac{s-\gamma/2}{\varepsilon^2}),$$
which gives the upper bound. \quad\quad $\blacksquare$
\vskip 0.4cm
\noindent{\bf Acknowledgement}. I thank the anonymous referee for the very helpful comments which have improved the exposition of the paper.


\begin{thebibliography}{99}
\bibitem{BDM}Budhiraja, A., Dupuis, P. and Maroulas, V. (2008). Large deviations for  infinite dimensional stochastic dynamical systems. {\it The Annals of Probability} {\bf 36}{4} 1390-1420.
\bibitem{CR}Cerrai, S., R\"{o}ckner,M. (2005). Large deviations for  invariant measures of stochastic reaction-diffusion systems with multiplicative noise and non-Lipschitz reaction term. {\it Ann. I. H. Poincar\'{e}.} {\bf 41} 69-105.
\bibitem{DMZ}Dalang, R., Mueller, C., Zambotti, L. (2006). Hitting properties of parabolic
SPDE's with reflection. {\it Ann. Probab.} {\bf 34}(4) 1423-1450.
\bibitem{DP1}Donati-Martin, C., Pardoux, E. (1993). White noise
driven SPDEs with reflection. {\it Probab. Theory Relat. Fields}
{\bf 95} 1-24.
\bibitem{DP2}Donati-Martin, C., Pardoux, E. (1997). EDPS r\'{e}fl\'{e}chies et calcul de
Malliavin. {\it Bull. Sci. Math.} {\bf 121}(5) 405-422.
\bibitem{FO}Funaki, T., Olla, S. (2001). Fluctuations for $\nabla\phi$ interface model on a wall.
 {\it Stoch. Proc. Appl.} {\bf 94}(1) 1-27.
 \bibitem{F}Freidlin, M. I., Wentzell, A. D. (1984). Random perturbations of dynamical systems.
 {\it Springer-Verlag} Berlin.
\bibitem{NP}Naulart, D., Pardoux, E. (1992). White noise driven quasilinear SPDEs with
reflection. {\it Probab. Theory Relat. Fields} {\bf 93} 77-89.
\bibitem{OS}Otobe, Y. (2006). Stochastic partial differential equations with two
reflecting walls. {\it J. Math. Sci. Univ. Tokyo} {\bf 13} 139-144.
\bibitem{S}Sowers, R. (1992). Large deviations for the invariant measure of a reaction-diffusion equation with non-Gaussian perturbations.  {\it Probability Theory and Related Fields} {\bf 92} 393-421.
\bibitem{W}Walsh, J. (1986). {\it An introduction to stochastic partial differential
equations.} Berlin Heidelberg Nerk York, Springer.
\bibitem{XZ}Xu, T., Zhang, T. (2009). White noise driven SPDEs with reflection: existence,
uniqueness and large deviation principles. {\it Stoch. Proc. Appl.} {\bf 119}(10) 3453-3470.
\bibitem{YZ1}Yang, J., Zhang, T. (2011). White noise driven SPDEs with two reflecting walls. {\it Infinite Dimensional Analysis, Quantum Probability and Related Topics.} {\bf 14}(4) 1-13.
\bibitem{YZ}Yang, J., Zhang, T. (2011). Existence and uniqueness of invariant measures of white noise driven SPDEs with two reflecting walls. Preprint 2011, Manchester.
\bibitem{ZA}Zambotti, L. (2001). A reflected stochastic heat equation as symmetric dynamic with respect
to the 3-d Bessel bridge.  {\it J. Funct. Anal.} {\bf 180} 195-209.
\bibitem{ZW}Zhang, T. (2010). White noise driven SPDEs with reflection: strong Feller properties and
 Harnack inequalities. {\it Potential Analysis.} {\bf 33:2} 137-151.
\end{thebibliography}
\end{document}